\documentclass[a4paper]{article}
\usepackage{graphicx}
\title{ Approximate Distribution of Hitting Probabilities 
for a Regular Surface with Compact Support in 2D}
\author{Denis S. Grebenkov}

\def\M{{\mathcal M}}
\def\E{{\mathcal E}}
\def\S{{\mathcal S}}
\def\iff{\textrm{if}\hskip 2mm }
\def\N{N^*}

\begin{document}

\begin{center}
\Large

{\bf  Approximate Distribution 
of Hitting Probabilities for 
a Regular Surface with Compact Support in 2D} 

\vskip 3mm
D.S.Grebenkov\footnote{ E-mail: Denis.Grebenkov@polytechnique.fr }
\normalsize
\vskip 3mm
Laboratoire de Physique de la Mati\`ere Condens\'ee, Ecole Polytechnique

91128 Palaiseau Cedex \hskip 1mm France
\vskip 2mm

Department of Statistical Physics, Saint Petersburg State University

ul. Ulyanovskaya 1, Petrodvorets, 198904, Saint Petersburg, Russia


\end{center}

\vskip 10mm

\begin{abstract}

Generalizing the well-known relations on characteristic
functions on a plane to the case of a one-dimensional regular 
surface (curve) with compact support,
we establish implicit equations for these functions.
Introducing an approximation, we solve the combinatorial problems
and reduce these equations to a set of linear 
equations for a finite number of unknown functions. 
Imposing natural conditions, we obtain a closed system of 
linear equations which can be solved for a given surface. 
Its solutions can be used to calculate the distribution of 
hitting probabilities for a regular surface with compact support.

In order to verify the accuracy of the approximate distribution
of hitting probabilities,
numerical analysis is being made for a chosen surface.

\end{abstract}

\vskip 10mm
\section*{ Introduction }

Search for the distribution of hitting probabilities is an 
old and a well-known problem. Consider random walk on $d$-dimensional 
lattice (in continuous case consider Brownian motion). 
Then fix a surface of interest $\S $. Suppose that any random walk 
starts from a given point $z$ which does not lie on $\S $. The
problem is to calculate the distribution $P_z(x)$ of probabilities 
of first contact with points $x$ of the surface $\S $. In other
words, we are looking for the probability that random walks 
from $z\notin \S$ to $x\in S$ do not touch other points 
$y\in \S \backslash \{ x\}$. Of course, the distribution $P_z(x)$
depends on $z$ and $\S $.

This problem has been solved exactly for some particular surfaces. 
For instance, the case of a planar surface in 2D (an ordinary
straight line) is described in any book on probability theory (see
\cite{Feller}, \cite{Walk}). Its generalization for $d$-dimensional
hyperplane is also simple (for example, see the end of Section 2). 
Note that exact solutions have been found only for some particular 
surfaces but not in the general case. In the general case, the asymptotic 
behavior is widely studied, \cite{Walk}.

Problems of the hitting probabilities do not only have a purely
mathematical interest. They are important for a wide class of
physical problems, in particular, for the problems of Laplacian 
transfer across an interface, for instance, diffusion through a 
membrane, electrod problems, heterogeneous catalysis, etc.
(see \cite{first}, \cite{1}, \cite{2} for details). Indeed,
if we are interested in diffusion through a semi-permeable
membrane (points of this membrane can absorb or reflect 
touching particle with certain probabilities), we can write
the total probability of absorption by a chosen point of the
membrane as a sum of probabilities to be absorbed 
after 0, 1, 2, etc. reflections (rigorous formalism
is described in \cite{2}, \cite{rapport}). Here we face
the task to calculate the distribution of hitting probabilities.
Note that using this distribution solely for a planar membrane, 
we have recently obtained some important results about general 
characteristics of the Laplacian transfer across an interface, 
\cite{my}.
To solve these problems one needs to know the distribution of 
hitting probabilities for a general surface. Here we propose a 
method to approximate the distribution of hitting probabilities 
for a rather general case in 2D.

In the first section we introduce definitions and conditions 
which are required in what follows. In the second section 
we briefly describe a well-known case of the hitting probabilities
on a horizontal axis. Main results are contained in the third section. 
Section 4 is devoted to some numerical results.
In the last section we make conclusions and discuss 
possible generalizations.

\vskip 10mm
\section{  Definitions }

Consider a square lattice on a plane.  Let us define a {\it regular 
surface\footnote{ Even for two-dimensional case we prefer to use the 
word ``surface'' instead of ``curve'' or something else.} with compact 
support} $\S =\{\: (x,S(x))\: \}$ by a function $S(x)$ with integer $x$ 
subject to the following conditions:

1. {\em Bijection}: The function $S(x)$ is a bijection between the set of integer
numbers (absiccae $x$) and the set of surface points ;

2. {\em Regularity}: For any $x$, \hskip 2mm $|S(x+1)-S(x)|\leq 1$ ;

3. {\em Compactness}: $\exists M: S(x)=0$ \hskip 2mm for $|x|\geq M$, i.e. 
the non-plane part of the surface has a finite size. In other words, function
$S(x)$ has a compact support. Moreover, we suppose that
the surface is centered : $S(\pm (M-1))\ne 0$.

\vskip 1mm

Let us briefly discuss this definition. 
The second condition allows to simplify all calculations
and formulae, but it does not seem to be essential (see Section 5). 
Note that this assumption can be viewed as {\it a regularity
condition} for the surface in continuous case : $S'(x)\leq 1$.

On the contrary, the third condition is important. It tells us
that the surface in question is a finite ``perturbation'' of a planar surface (line).
In other words, this surface is composed of two parts: 
a complex but compact part in the center with two plane ``tails''.
Moreover, it is important that both tails lie on the same height
(which is chosen as $0$). This feature will allow to obtain 
an approximate distribution of hitting probabilities by using 
the same ideas as for a planar surface (see Section 2). 

\vskip 1mm

We call all the points $\{ (x,y)\; :\; y=n\}$ the {\it $n^{th}$ level}. 
Denote 
$$N=\max \{ S(x) \}, \hskip 15mm \N=-\min \{ S(x) \} ,$$ 
i.e. the surface lies between $(-\N)$th and $N$th levels.

All points $\M=\{ (x,y)\; :\; \forall x\hskip 3mm y<S(x)\}$ are called
{\it internal}. All points $\E=\{ (x,y)\; :\; \forall x\hskip 3mm y>S(x)\}$ 
are called {\it external}. The external points near the surface, $\{ (x,S(x)+1) \}$ 
are called {\it near-boundary points}. The functions defined on these
points, are called {\it near-boundary functions} (see below).
Often we shall use the words ``surface'', ``near-boundary functions'', etc.
thinking only about the non-trivial part, i.e. for $|x|<M$.

The external points with $y=0$ are called 
{\it ground points}. The functions defined on these points, 
are called {\it ground functions}. Let $J=\{ k\: : \: (k,0)\in \E \}$
the set of abscissae of ground points. Let also 
$J_0=\{ k\in (-M,M)\: : \: (k,0)\in S \}$ the set of abscissae 
of boundary points on zeroth level (only non-plane part!).

\vskip 1mm

We introduce the hitting probabilities $P_{k,n}(x)$, i.e.
the probability of the first contact with the surface at point $(x,S(x))$ if started 
from $(k,n)$. Their characteristic functions $\phi _{m,n}(\theta )$ are
\begin{equation}    \label{phi_def}
\phi _{k,n}(\theta )=\sum\limits _{x=-\infty }^{\infty }P_{k,n}(x)
e^{ix\theta }  .
\end{equation}
The inverse Fourier transform allows to obtain $P_{k,n}(x)$,
\begin{equation}   \label{Fourier}
P_{k,n}(x)=\int\limits _{-\pi }^{\pi }\frac{d\theta }{2\pi }e^{-ix\theta }
\phi _{k,n}(\theta ) .
\end{equation}

\vskip 10mm
\section{ Planar surface }

At the beginning, we consider the trivial and well-known case of a 
planar surface (horizontal axis): $S(x)=0$. This case is useful to
remind the technique of manipulation with characteristic functions.

Suppose that $n>0$. The probability $P_{k,n}(x)$ satisfies a simple identity
\begin{equation}   \label{Pidentity}
P_{k,n}(x)=\frac14\biggl[P_{k+1,n}(x)+P_{k-1,n}(x)+P_{k,n+1}(x)+P_{k,n-1}(x)\biggr] ,
\end{equation}
which can be also written for characteristic functions,
\begin{equation}    \label{identity}
\phi _{k,n}(\theta )=\frac14\biggl[\phi _{k+1,n}(\theta )+\phi _{k-1,n}(\theta )+
\phi _{k,n+1}(\theta )+\phi _{k,n-1}(\theta )\biggr] .
\end{equation}
Translational invariance along the horizontal axis gives
\begin{equation}   \label{trans}
\phi _{k,n}(\theta )=e^{ik\theta }\phi _{0,n}(\theta ).
\end{equation}
Using the obvious condition $P_{k,0}(x)=\delta _{k,x}$, we obtain
\begin{equation}   \label{phi0_0}
\phi _{k,0}(\theta )=e^{ik\theta } .
\end{equation}
The last trick is the following. If the starting point is placed in the
$n$-th level, the random walk must cross the $(n-1)$-th level at some point 
$(m,n-1)$ to reach zeroth level. The probability to pass from $(k,n)$ to 
$(m,n-1)$ without touching other points in the $(n-1)$-th level is exactly 
$P_{k,1}(m)$. Therefore we can write
$$P_{k,n}(x)=\sum\limits _m P_{k,1}(m)P_{m,n-1}(x) .$$
In terms of characteristic functions this convolution is
just a product of the two corresponding characteristic functions,
$$\phi _{k,n}(\theta )=\phi _{0,1}(\theta )\phi _{k,n-1}(\theta ).$$
Using the translational invariance (\ref{trans}), we obtain
\begin{equation}    \label{phi0_n}
\phi _{k,n}(\theta )=e^{ik\theta }[\phi _{0,1}(\theta )]^n .
\end{equation}
Substitution of expressions (\ref{phi0_0}) and (\ref{phi0_n}) into relation 
(\ref{identity}) for $n=1$ and $k=0$ leads to
$$\phi _{0,1}(\theta )=\frac14\biggl(e^{-i\theta }\phi _{0,1}(\theta )+
e^{i\theta }\phi _{0,1}(\theta )+1+[\phi _{0,1}(\theta )]^2\biggr) ,
\hskip 5mm \textrm{or}$$
\begin{equation}   \label{eqn_plane}
\phi _{0,1}^2-(4-2\cos \theta )\phi _{0,1}+1=0.
\end{equation}
This quadratic equation has two solutions, and we should choose the 
one for which $\phi _{0,1}(\theta )\leq 1$ (property of characteristic 
function). It is denoted $\varphi (\theta )$,
\begin{equation}    \label{varphi}
\varphi (\theta )=2-\cos \theta -\sqrt{(2-\cos \theta )^2-1} .
\end{equation}
So, we obtain for the planar surface
\begin{equation}   \label{phi_plane}
\phi _{k,n}(\theta )=e^{ik\theta }\varphi ^n(\theta ).
\end{equation}
Inverting this relation with the help of (\ref{Fourier}), we obtain
the distribution of hitting probabilities for the planar surface,
\begin{equation}   \label{H}
P^{planar}_{k,n}(x)=\int\limits _{-\pi }^{\pi }\frac{d\theta }{2\pi }
e^{i(k-x)\theta }\varphi ^n(\theta )=H^n_{k-x} .
\end{equation}
This well-known result will be used for a general case (some
properties of coefficients $H^n_k$ are described in Appendix 1).
The formulae (\ref{varphi}) and (\ref{phi_plane}) can be generalized
for $d$-dimentional hypercubic lattice,
$$\varphi (\theta _1...\theta _{d-1})=d-\sum\limits _{i=1}^{d-1}\cos (\theta _i)
-\sqrt{\left(d-\sum\limits _{i=1}^{d-1}\cos (\theta _i)\right)^2-1} ,$$
$$\phi _{{\bf x},n}(\theta _1...\theta _{d-1})=\exp \biggl[i\sum\limits _{i=1}^{d-1}
x_i\theta _i]\biggr]\varphi ^n(\theta _1...\theta _{d-1}) .$$

\vskip 10mm
\section{ Regular surface with compact support }

We shall consider the characteristic functions $\phi _{k,n}(\theta )$
as a vector 
$$\Phi ^{(n)}(\theta )=\left(\begin{array}{ c }
\phi _{-L,n}(\theta )  \\
\phi _{-L+1,n}(\theta ) \\
  ...    \\
\phi _{L,n}(\theta ) \\
\end{array}\right)$$
of $(2L+1)$ components where parameter $L$ is supposed large,
and it will tend to infinity at the end of calculation.

For the planar surface we had relation (\ref{identity}) which
can be written in matrix form
\begin{equation}    \label{A0}
A\Phi ^{(n)}=\Phi ^{(n-1)}+\Phi ^{(n+1)},
\end{equation}
where the matrix $A$ is
$$A_{i,i}=4, \hskip 5mm  A_{i,i+1}=A_{i+1,i}=-1,  
\hskip 5mm   A_{2L+1,1}=A_{1,2L+1}=-1  .$$
The last equalities are artificial : we added them to obtain a 
cyclic structure of $A$. But at the limit $L\to \infty $ this 
little modification vanishes. The eigenvalues of $A$ are 
$$\lambda _h=4-2\cos (\theta _h), \hskip 5mm \textrm{with}   \hskip 2mm
\theta _h=\frac{2\pi h}{2L+1} ,$$
and the eigenvectors are given as
$$V_h=\left(\begin{array}{ c }
e^{-iL\theta _h }  \\
e^{-i(L-1)\theta _h } \\
  ...    \\
e^{iL\theta _h } \\
\end{array}\right).$$
Now we generalize the relation (\ref{A0}) to the case of
a regular surface with compact support by introducing vector 
$\Delta \Phi ^{(n)}$,
\begin{equation}    \label{A}
A\Phi ^{(n)}=\Phi ^{(n-1)}+\Phi ^{(n+1)}+\Delta \Phi ^{(n)}
\end{equation}
(this relation can be regarded as the definition of $\Delta \Phi ^{(n)}$).
Let us introduce 
\begin{equation}    \label{c_mn}
c_n(\theta ,\theta _h)=(\Phi ^{(n)},V^*_h), \hskip 10mm  
\Delta c_n(\theta ,\theta _h)=(\Delta \Phi ^{(n)},V^*_h) . 
\end{equation}
We can rewrite (\ref{A}) in terms of $c_n$ and $\Delta c_n$,
\begin{equation}    \label{spect}
\lambda _hc_n=c_{n-1}+c_{n+1}+\Delta c_n .
\end{equation}
If we can express $c_n$ in terms of $\lambda _h$, $\varphi (\theta )$,
$c_0$ and $\{ \Delta c_k\}$, we find $\Phi ^{(n)}$ as a decomposition 
in the eigenbasis $V_h$,
\begin{equation}    \label{solution}
\Phi ^{(n)}=\frac{1}{2L+1}\sum\limits _h c_n(\theta ,\theta _h)V_h  
\end{equation}
(factor $(2L+1)^{-1}$ is due to normalization $(V_h,V^*_h)$). 
In order to solve the recurrence relations (\ref{spect}), we should 
close them by certain conditions. For the lower half plane we take 
a sufficiently large number $N_l>\N $, and
\begin{equation}    \label{N*}
\Phi _{-N_l}=0 , \hskip 10mm \textrm{or}  \hskip 10mm  c_{-N_l}=0 ,
\end{equation}
because it is impossible to penetrate through the surface.
More generally, according to the definition of hitting probabilities, 
we should maintain
\begin{equation}      \label{zeros}
\phi _{m,n}(\theta )=\cases{ e^{im\theta }, \: \iff (m,n)\in \S ,\cr
                             0, \hskip 6mm  \iff (m,n)\in \M . }
\end{equation}
For the upper half plane we shall use the following trick. 
We take a large number $N_u\gg N$ and consider $\phi _{x,N_u+1}(\theta )$.
As for a planar surface, we can write $P_{x,N_u+1}(n)$ in convolution form,
$$P_{x,N_u+1}(n)=\sum\limits _m H^1_{x-m}P_{m,N_u}(n) ,$$
where $H^1_{x-m}$ is the probability to hit point $(m,N_u)$ if started from 
$(x,N_u+1)$. For characteristic functions it is simply
$$\phi _{x,N_u+1}(\theta )=\sum\limits _m H^1_{x-m}\phi _{m,N_u}(\theta ) .$$
In the plane case we used the translational invariance along the horizontal axis
to simplify this sum. Evidently, such a symmetry breaks down in the 
general case. But if we take $N_u$ sufficiently large, i.e. we ``look''
on the membrane from a remoted point, we can suppose that translational
invariance is approximately valid,
\begin{equation}   \label{real_approx}
\phi _{m,N_u}(\theta )\approx e^{i(m-x)\theta }\phi _{x,N_u}(\theta ) .
\end{equation}
Using this approximation, we immediately obtain 
$$\phi _{x,N_u+1}(\theta )=\varphi (\theta )\phi _{x,N_u}(\theta ) ,$$
where function $\varphi $ is defined by (\ref{varphi}).
This is our approximation in the upper half plane which allows to
close the recurrence relations (\ref{spect}),
\begin{equation}     \label{N+1}
c_{N_u+1}(\theta ,\theta _h)=\varphi (\theta )c_{N_u}(\theta ,\theta _h) .
\end{equation}

The main idea is to step down from $N_u$-th and $(-N_l)$-th levels to zeroth level.
We shall consider the upper and lower half planes separately because the relations
(\ref{N*}) and (\ref{N+1}) are different.
Note the essential complication of the general case with respect to
the planar surface. For the planar surface we had $\Delta c_n=0$ for any $n$, 
and the system of equations (\ref{spect}) was closed. It was sufficient
to solve these recurrence relations by substitution $c_n=c_0c^n$ in
(\ref{spect}), and we obtained the final form of $\Phi ^{(n)}$.
On the contrary, for the general surface $\Delta c_n\ne 0$,
and they depend on the near-boundary functions $\phi _{m,n}$ (see below). 
Consequently, the decomposition (\ref{solution}) itself becomes a system of 
implicit equations for $\phi _{m,n}$. For the moment, the  problem is 
complex. It will be solved by two steps. First, we obviate the combinatorial 
problems, i.e. we find the explicit solution of reccurence relations (\ref{spect}).
Second, we solve the equations for $\phi _{m,n}$.

\vskip 5mm 
\subsection{ Solution of reccurence relations. }

A direct verification shows that
\begin{eqnarray}   \label{rec}
c_n &=& \beta _{N_u-n}c_{N_u}-\sum\limits _{l=1}^{N_u-n}\alpha _l\Delta c_{n+l} , \\  \label{rec*}
c_{-n} &=& \alpha _{N_l-n}c_{-N_l+1}-\sum\limits _{l=1}^{N_l-n-1}\alpha _l\Delta c_{-n-l} 
\end{eqnarray}
is a general solution of (\ref{spect}) (we omitted the index $h$ 
which does not change the structure of the solution), where
\begin{equation}   \label{alpha}
\alpha _0=0, \hskip 4mm  \alpha _1=1, \hskip 10mm  
\alpha _{n+2}(\theta _h)=\lambda _h\alpha _{n+1}(\theta _h)-\alpha _n(\theta _h) ,
\end{equation}
\begin{equation}   \label{beta}
\beta _n(\theta ,\theta _h)=\alpha _n(\theta _h)\bigl[\lambda _h-\varphi (\theta )\bigr]-
\alpha _{n-1}(\theta _h) ,  \hskip 6mm  \beta _0=1 .
\end{equation}
We are looking for the explicit representation for $\alpha _n$ in the form 
$\alpha _n=x^{n-1}$. Substituting this into (\ref{alpha}), we obtain equation
$$x^2-\lambda _h x+1=0, $$
which has two well-known solutions: $\varphi $ and $\varphi ^{-1}$
(compare this equation with (\ref{eqn_plane})). As the expression
(\ref{alpha}) is linear, we find a general solution as linear
combination of $\varphi ^{n-1}$ and $\varphi ^{1-n}$ such that
$\alpha _0=0$. We obtain
\begin{equation}   \label{alpha1}
\alpha _n(\theta _h)=\frac{1-\varphi ^{2n}(\theta _h)}{1-\varphi ^2(\theta _h)}
\varphi ^{1-n}(\theta _h) ,
\end{equation}
or as a geometrical sequence
\begin{equation}   \label{alpha2}
\alpha _n(\theta _h)=\sum\limits _{j=0}^{n-1}\varphi ^{2j+1-n}(\theta _h) .
\end{equation}

Formulae (\ref{rec}) and (\ref{rec*}) are valid for any $n\geq 0$, in particular, 
for $n=0$, and we can express $c_{N_u}$ and $c_{-N_l+1}$ in terms of $c_0$ and 
$\{ \Delta c_l\}$,
$$c_{N_u}=\frac{1}{\beta _{N_u}}\left(c_0+\sum\limits _{l=1}^{N_u}\alpha _l\Delta c_l\right) ,\hskip 10mm
c_{-N_l+1}=\frac{1}{\alpha _{N_l}}\left(c_0+\sum\limits _{l=1}^{N_l-1}\alpha _l\Delta c_{-l}\right) ,$$
hence
\begin{eqnarray}   \label{solu}
c_n &=& \frac{\beta _{N_u-n}}{\beta _{N_u}}\left(c_0+\sum\limits _{l=1}^{N_u}\alpha _l\Delta c_l\right)
-\sum\limits _{l=1}^{N_u-n}\alpha _l\Delta c_{n+l} , \\   \label{solu*}
c_{-n} &=& \frac{\alpha _{N_l-n}}{\alpha _{N_l}}\left(c_0+\sum\limits _{l=1}^{N_l-1}
\alpha _l\Delta c_{-l}\right)-\sum\limits _{l=1}^{N_l-n-1}\alpha _l\Delta c_{-n-l} .
\end{eqnarray}
Let us introduce 
$$f^{(N_u)}_n(\theta ,\theta _h)=\frac{\beta _{N_u-n}}{\beta _{N_u}}, \hskip 15mm
\tilde{f}^{(N_l)}_n(\theta _h)=\frac{\alpha _{N_l-n}}{\alpha _{N_l}}  $$
(the right-hand side of $f^{(N_u)}_n$ depends on $\theta $ through factor $\varphi $;
the dependence on $\theta _h$ is due to $\alpha _n$ which contains $\lambda _h$).

Using the reccurence properties of $\alpha _n$ (see Appendix 2), 
we simplify relations (\ref{solu}) and (\ref{solu*}),
\begin{eqnarray}    \label{cn}
c_n &=& f^{(N_u)}_n\left(c_0+\sum\limits _{l=1}^n\alpha _l\Delta c_l\right)+\alpha _n
\sum\limits _{l=n+1}^{N_u}f^{(N_u)}_l\Delta c_l , \\   \label{cn*}
c_{-n} &=& \tilde{f}^{(N_l)}_n\left(c_0+\sum\limits _{l=1}^n\alpha _l\Delta c_{-l}\right)+
\alpha _n\sum\limits _{l=n+1}^{N_l-1}\tilde{f}^{(N_l)}_l\Delta c_{-l} . 
\end{eqnarray}

Rewriting the definition (\ref{beta}) of $\beta _n$ as
$$\beta _n(\theta, \theta _h)=\alpha _{n+1}(\theta _h)-
\varphi (\theta )\alpha _n(\theta _h) ,$$
after simplifications we obtain
$$f^{(N_u)}_n(\theta ,\theta _h)=\frac{\varphi ^n(\theta _h)[1-\varphi (\theta )
\varphi (\theta _h)]-\varphi ^{2N_u-n+1}(\theta _h)[\varphi (\theta _h)-
\varphi (\theta )]}{[1-\varphi (\theta )\varphi (\theta _h)]
-\varphi ^{2N_u+1}(\theta _h)[\varphi (\theta _h)-\varphi (\theta )]} .$$
We remind that $N_u$ is an {\it arbitrary} sufficiently large number. 
Therefore we can take the limit $N_u\to \infty $. Knowing that 
$\varphi (\theta _h)<1$ for $\theta _h\ne 0$, we obtain that in this limit 
$\varphi ^{2N_u-n+1}$ and $\varphi ^{2N_u+1}$ vanish, i.e.
\begin{equation}    \label{approx}
f^{(\infty )}_n(\theta ,\theta _h)=\varphi ^n(\theta _h) . 
\end{equation}
Using formula (\ref{alpha1}), we can also write the explicit 
representation for $\tilde{f}^{(N_l)}_n$,
$$\tilde{f}^{(N_l)}_n(\theta _h)=\frac{1-\varphi ^{2(N_l-n)}(\theta _h)}
{1-\varphi ^{2N_l}(\theta _h)}\varphi ^n(\theta _h) .$$
As above, we take the limit $N_l\to \infty $ to obtain
\begin{equation}    \label{approx*}
\tilde{f}^{(\infty )}_n(\theta _h)=\varphi ^n(\theta _h)  .
\end{equation}
Note that we obtained the same limits $f^{(\infty )}_n$ and $\tilde{f}^{(\infty )}_n$
imposing the different conditions (\ref{N*}) and (\ref{N+1}) in the lower and
upper half planes respectively. Using formulae (\ref{approx}) and (\ref{approx*}),
we can write
\begin{eqnarray*}    
c_n &=& \varphi ^n(\theta _h)\left(c_0(\theta, \theta _h)+\sum\limits _{l=1}^n
\alpha _l(\theta _h)\Delta c_l(\theta ,\theta _h)\right)+\alpha _n(\theta _h)
\sum\limits _{l=n+1}^{\infty }\varphi ^l(\theta _h)\Delta c_l(\theta ,\theta _h) , \\   
c_{-n} &=& \varphi ^n(\theta _h)\left(c_0(\theta ,\theta _h)+\sum\limits _{l=1}^n\alpha _l(\theta _h)
\Delta c_{-l}(\theta ,\theta _h)\right)+\alpha _n(\theta _h)\sum\limits _{l=n+1}^{\infty }
\varphi ^l(\theta _h)\Delta c_{-l}(\theta ,\theta _h) . 
\end{eqnarray*}
Introducing functions
\begin{equation}    \label{gamma}
\gamma ^{(n)}_l(\theta )=\cases{ 
   \varphi ^n(\theta )\alpha _l(\theta ),\hskip 2mm \iff l\leq n, \cr
   \varphi ^l(\theta )\alpha _n(\theta ),\hskip 2mm \iff l>n , \cr 
  \hskip 8mm    0, \hskip 9mm  \iff l\leq 0 \hskip 2mm \textrm{or}  \hskip 2mm  n\leq 0}
\end{equation}
(the last convention will be used in the following), we can 
write $c_n$ in the unique form (for $n>0$ and $n<0$),
\begin{equation}    \label{cn3}
c_n=\varphi ^{|n|}(\theta _h)c_0(\theta, \theta _h)+\sum\limits _{l=-\N}^{N}
\bigl[\gamma ^{(n)}_l(\theta _h)+\gamma ^{(-n)}_{-l}(\theta _h)\bigr]
\Delta c_l(\theta ,\theta _h) ,
\end{equation}
where summation over $l$ is just from $-\N $ to $N$, because corrections
$\Delta c_{N+l}$ and $\Delta c_{-\N-l}$ are equal to $0$ for $l>0$ 
(see Section 3.4 for details).

\vskip 5mm
\subsection{ Coefficient $c_0$ }
  
Let us calculate the coefficient $c_0$,
$$c_0(\theta ,\theta _h)=(\Phi _0,V^*_h)=\sum\limits _{k=-L}^Le^{-ik\theta _h}\phi _{k,0}
=c^{(0)}_0+c^{(1)}_0.$$
As it was mentioned above, the plane ``tails'' of the surface
lie on the zeroth level, thus
$$(\Phi _0)_k=e^{ik\theta },\hskip 2mm \iff |k|\geq M.$$
Using this explicit form, we are going to compute the contribution $c_0^{(0)}$
of plane ``tails'' 
\begin{equation}        \label{c0}
c^{(0)}_0(\theta ,\theta _h)=2\sum\limits _{k=M}^L\cos (\theta -\theta _h)=
\frac{\sin (L+0.5)(\theta -\theta _h)}{\sin (\theta -\theta _h)/2}-
\frac{\sin (M-0.5)(\theta -\theta _h)}{\sin (\theta -\theta _h)/2}
\end{equation}
(in the case $\theta =\theta _h$ one should consider this relation
in the limit sense when $\theta \to \theta _h$).
Note the simple relation
\begin{equation}      \label{integral}
\int\limits _{-\pi /2}^{\pi /2}\frac{d\theta }{\pi }
\cos (2k\theta )\frac{\sin (2M-1)\theta }{\sin \theta }=\chi _{(-M,M)}(k) ,
\hskip 10mm     \chi _A(k)=\cases{  1, \hskip 2mm \iff k\in A , \cr
                      0, \hskip 2mm \iff k\notin A .}
\end{equation}
%

The contribution $c_0^{(1)}$ of intermediate 
part of the surface (with $|x|<M$) can contain some non-trivial
functions $\phi _{m,0}(\theta )$, with $k\in J$ and $k\in J_0$
\begin{equation}    \label{c01}
c_0^{(1)}(\theta ,\theta _h)=\sum\limits _{k\in J}\phi _{k,0}(\theta )e^{-ik\theta _h}
+\sum\limits _{k\in J_0}e^{ik(\theta -\theta _h)} .
\end{equation}

\vskip 5mm
\subsection{ Limit $L\to \infty $ }

Expression (\ref{cn3}) transforms to the integral expression 
for $\phi _{k,n}$ by taking the limit $L\to \infty $. Here we write 
only the expression for $n\geq 0$, the opposite case will be easily 
obtained later. Note that the first term of (\ref{c0})
tends to $\delta $-function in the limit $L\to \infty $. It removes 
the integration in the first term, i.e. we obtain
\begin{equation}    \label{solu1}
\phi _{k,n}(\theta )=\varphi ^n(\theta )e^{ik\theta }-
\phi ^{(0)}_{k,n}(\theta )+
\end{equation}
$$+\int\limits _{-\pi }^\pi \frac{d\theta '}{2\pi }e^{ik\theta '}\left(
\varphi ^n(\theta ')c_0^{(1)}(\theta ,\theta ')+
\sum\limits _{l=1}^{N}\gamma ^{(n)}_l(\theta ')\Delta c_l(\theta ,\theta ')\right) .$$
The first term is the contribution of plane ``tails''. The second term
corresponds to the perturbation on zeroth level 
due to $c_0^{(0)}$, 
\begin{equation}    \label{phi0}
\phi ^{(0)}_{k,n}(\theta )=\int\limits _{-\pi }^\pi \frac{d\theta '}{2\pi }
e^{ik\theta '}\varphi ^n(\theta ')\frac{\sin (M-0.5)(\theta -\theta ')}
{\sin (\theta -\theta ')/2} .
\end{equation}
The last term of (\ref{solu1}) contains some unknown functions 
$\phi _{m,l}(\theta )$ through the coefficients $\Delta c_l(\theta ,\theta ')$. 
It is denoted as $T[\phi ]$, and we are going to calculate it.

\vskip 5mm
\subsection{ Coefficients $\Delta c_l$ }

To get ahead with the expression (\ref{solu1}), we should write
explicitly the coefficients $\Delta c_l(\theta ,\theta ')$. It is not so 
easy for the general case. Indeed, for these purposes one  
calculates the contributions of each point on $l$-th level. 
The problem is that there are many conditions, and they lead to complex 
formulae difficult to manipulate. We are going to present the other way.

What is the origin of the vector $\Delta \Phi ^{(n)}$?  Let us recall 
the definition (\ref{A}) where vectors $\Delta \Phi ^{(n)}$ were introduced 
to generalize the expression (\ref{A0}). A brief reflection shows that
$$(\Delta \Phi ^{(n)})_m=\cases{  0, \hskip 70mm \iff  (m,n)\in \E , \cr
  4\phi _{m,n}-\phi _{m-1,n}-\phi _{m+1,n}-\phi _{m,n-1}-\phi _{m,n+1}, 
         \hskip 2mm   \iff  (m,n)\notin \E }.$$
In other words, the relation (\ref{A0}) is satisfied automatically
for any external point, but it should be imposed artificially for each 
surface and internal points.

Now it is the moment to remind the formula (\ref{zeros}) which
tells that functions $\phi _{m,n}(\theta )$ are equal to zero
on the internal points. Therefore we can consider only the points
near the surface $\S $. A direct verification shows that
$$(\Delta \Phi ^{(S(m))})_m=4e^{im\theta }-\phi _{m,S(m)+1}-\cases{
  e^{i(m+1)\theta }, \hskip 5mm   \iff  S(m+1)-S(m)=0 \cr
  \phi _{m+1,S(m)},  \hskip 2mm   \iff  S(m+1)-S(m)=-1 \cr
  0,     \hskip 16mm              \iff  S(m+1)-S(m)=1 } $$
\begin{equation}   \label{In1}
-\cases{
  e^{i(m-1)\theta }, \hskip 5mm   \iff  S(m-1)-S(m)=0  \cr
  \phi _{m-1,S(m)},  \hskip 2mm   \iff  S(m-1)-S(m)=-1  \cr
  0,    \hskip 16mm               \iff  S(m-1)-S(m)=1 } ,
\end{equation}
$$(\Delta \Phi ^{(S(m)-1)})_m=-e^{im\theta }-e^{i(m+1)\theta }\delta _{S(m),S(m+1)+1}
-e^{i(m-1)\theta }\delta _{S(m),S(m-1)+1} , $$
$$(\Delta \Phi ^{(n)})_m=0, \hskip 5mm \iff n\ne S(m) 
\hskip 2mm \textrm{and} \hskip 2mm n\ne S(m)-1 $$
(here we use Kronecker $\delta $-symbols, $\delta _{ii}=1$ and $\delta _{ij}=0$ if $i\ne j$).
Usually there are several nonzero components of $\Delta \Phi ^{(n)}$ for each $n$,
because each level contains several surface points. But there exists
one exception -- zeroth level, where there is infinity of 
surface points due to the plane ``tails''. Thus, the vector 
$\Delta \Phi ^{(-1)}$ has exceptional structure\footnote{ There exists another 
vector with exceptional structure, $\Delta \Phi ^{(0)}$. However, we construct 
our treatment by the way when this vector does not appear in expressions (see (\ref{cn}), 
(\ref{cn*})). In other words, we step down from $N_u$-th and $(-N_l)$-th levels to zeroth level
in order to do not pass through the zeroth level.}.
It contains the usual terms due to the non-trivial part of the surface, 
and the contribution of plane ``tails''. Note that the last one is equal 
to $-c_0^{(0)}$ which was calculated in Section 3.2. Later we shall use
this result for the lower half plane.

\vskip 5mm
\subsection{  Approximate distribution of hitting probabilities }

According to the definition (\ref{c_mn}) of $\Delta c_l$, we can rewrite
$T[\phi ]$ as
$$T[\phi ]=\int\limits _{-\pi }^\pi \frac{d\theta '}{2\pi }e^{ik\theta '}
\sum\limits _{m=-M+1}^{M-1}e^{-im\theta '}\sum\limits _{l=1}^{N}
\gamma ^{(n)}_l(\theta ')(\Delta \Phi ^{(l)})_m ,$$
where we changed the order of summation over $m$ and $l$.
However, in the last sum there are only two terms corresponding to
$(\Delta \Phi ^{(S(m))})_m$ and $(\Delta \Phi ^{(S(m)-1)})_m$, if $S(m)\geq 1$ 
(otherwise, this sum is equal to $0$). Using expression (\ref{In1}), we obtain 
explicitly
$$T[\phi ]=\sum\limits _{m=-M+1}^{M-1}\int\limits _{-\pi }^\pi \frac{d\theta '}{2\pi }
e^{i(k-m)\theta '}\biggl(\gamma ^{(n)}_{S(m)}(\theta ')
\biggl[4e^{im\theta }-\phi _{m,S(m)+1}-e^{i(m+1)\theta }\delta _{S(m),S(m+1)}$$
$$-\phi _{m+1,S(m)}\delta _{S(m),S(m+1)+1}-e^{i(m-1)\theta }\delta _{S(m),S(m-1)}
-\phi _{m-1,S(m)}\delta _{S(m),S(m-1)+1}\biggr]$$
$$-\gamma ^{(n)}_{S(m)-1}(\theta ')
\biggl[e^{im\theta }+e^{i(m-1)\theta }\delta _{S(m),S(m-1)+1}+
e^{i(m+1)\theta }\delta _{S(m),S(m+1)+1}\biggr]\biggr).$$
(here we have used the last convention in the definition (\ref{gamma}) of 
$\gamma ^{(n)}_l$ to avoid any terms with $S(m)\leq 0$).

The last step is to transform this huge expression for characteristic
functions into hitting probabilities using the formula (\ref{Fourier}).
Note that all functions $e^{im\theta }$ after integration over 
$\theta $ with $e^{-ix\theta }$ give $\delta $-symbols that remove the
summation over $m$ in corresponding terms (but we should 
write factor $\chi _{(-M,M)}(x)$),
\begin{equation}    \label{supplem}
T[P]=-\sum\limits _{m=-M+1}^{M-1}D^{(k,n)}_{m,S(m)}\biggl(P_{m,S(m)+1}(x)+
\end{equation}
$$P_{m+1,S(m)}(x)\delta _{S(m),S(m+1)+1}
+P_{m-1,S(m)}(x)\delta _{S(m),S(m-1)+1}\biggr)+$$
$$\chi _{(-M,M)}(x)\biggl(4D^{(k,n)}_{x,S(x)}-D^{(k,n)}_{x-1,S(x-1)}
\delta _{S(x),S(x-1)}-D^{(k,n)}_{x+1,S(x+1)}\delta _{S(x),S(x+1)}-$$
$$D^{(k,n)}_{x,S(x)-1}-D^{(k,n)}_{x+1,S(x+1)-1}\delta _{S(x),S(x+1)-1}
-D^{(k,n)}_{x-1,S(x-1)-1}\delta _{S(x),S(x-1)-1}\biggr) ,$$
where 
\begin{equation}    \label{D0}
D^{(k,n)}_{m,l}=\int\limits _{-\pi }^\pi \frac{d\theta '}{2\pi }
e^{i(k-m)\theta '}\gamma ^{(n)}_l(\theta ') .
\end{equation}
Note that coefficients $D^{(k,n)}_{m,l}$ are universal, they
do not depend on a given surface. It means that once calculated,
these coefficients can be used for any hitting problem in 2D.
They can be also expressed in terms of $H^n_k$,
$$D^{(k,n)}_{m,l}=\sum\limits _{j=1}^{\min \{n,l\}}H^{2j-1+|l-n|}_{k-m}$$ 
(if $n$ or $l$ is equal to $0$, the sum is also equal to $0$).
We just indicate several useful properties of these coefficients,
$$D^{(k,n)}_{m,l}=D^{(0,n)}_{m-k,l}=D^{(k-m,n)}_{0,l}=D^{(m-k,n)}_{0,l},
\hskip 10mm D^{(k,0)}_{m,l}=D^{(k,n)}_{m,0}=0 ,$$
$$D^{(k,n)}_{m+1,l}+D^{(k,n)}_{m-1,l}+D^{(k,n)}_{m,l+1}+D^{(k,n)}_{m,l-1}=
4D^{(k,n)}_{m,l}-\delta _{k,m}\delta _{n,l} .$$
  
\vskip 1mm

The first part of (\ref{supplem}), containing $P_{m,S(m)+1}$, can be represented as
$$T_1[P]=-\sum\limits _{m=-M}^{M}G^{(k,n)}_mP_{m,S(m)+1}(x)$$
with coefficients 
$$G^{(k,n)}_m=D^{(k,n)}_{m,S(m)}+\delta _{S(m),S(m+1)-1}D^{(k,n)}_{m+1,S(m+1)}
+\delta _{S(m),S(m-1)-1}D^{(k,n)}_{m-1,S(m-1)} .$$

The second part of (\ref{supplem}) can be simplified. Indeed,
using the properties of $\delta $-symbols, we have
$$T_2=\chi _{(-M,M)}(x)\biggl(4D^{(k,n)}_{x,S(x)}-D^{(k,n)}_{x-1,S(x)}-
D^{(k,n)}_{x+1,S(x)}-D^{(k,n)}_{x,S(x)-1}$$
$$+D^{(k,n)}_{x-1,S(x)}\delta _{S(x),S(x-1)+1}+D^{(k,n)}_{x+1,S(x)}
\delta _{S(x),S(x+1)+1}\biggr) .$$
Using the properties of $D^{(k,n)}_{m,l}$, we finally obtain
\begin{equation}    \label{T2}
T_2=\chi _{(-M,M)}(x)\biggl(D^{(k,n)}_{x,S(x)+1}
+D^{(k,n)}_{x-1,S(x)}\delta _{S(x),S(x-1)+1}+D^{(k,n)}_{x+1,S(x)}
\delta _{S(x),S(x+1)+1}\biggr) .
\end{equation}
Here we omitted term $\delta _{k,x}\delta _{n,S(x)}$ supposing that
starting point $(k,n)$ does not lie on the surface.

\vskip 1mm

Let us get back to the formula (\ref{solu1}). Using the inverse 
Fourier transform (\ref{Fourier}), we write
\begin{equation}   \label{prob}
P_{k,n}(x)=H^n_{k-x}-P^{(0)}_{k,n}(x)+\sum\limits _{m\in J_0}H^n_{k-m}\delta _{m,x}+
\sum\limits _{m\in J}H^n_{k-m}P_{m,0}(x)+T_1[P]+T_2  .
\end{equation}
The second term is 
$$P_{k,n}^{(0)}(x)=\int\limits _{-\pi }^{\pi }\frac{d\theta }{2\pi }e^{-ix\theta }
\int\limits _{-\pi }^{\pi }\frac{d\theta '}{2\pi }e^{ik\theta '}
\varphi ^n(\theta ')\frac{\sin (M-0.5)(\theta -\theta ')}
{\sin (\theta -\theta ')/2} .$$
Replacing in the first integral $\theta _1=\theta -\theta '$, we
factorize these integrals.
The first factor is exactly $H^n_{k-x}$. The second one was calculated
explicitly, see (\ref{integral}), and it equals to $\chi _{[-M+1,M-1]}(x)$.
Consequently, the first two terms in (\ref{prob}) can be grouped into
$H^n_{k-x}\chi _{(-\infty,-M]\cap [M,+\infty )}(x)$. It means that 
the solution $H^n_{k-x}$ of the plane case is valid only for the plane 
``tails'', whereas on the non-trivial surface (for $|x|<M$) the main 
contribution is due to other terms.  So, we have obtained an important result,
\begin{equation}   \label{prob2}
P_{k,n}(x)=\tilde{P}_{k,n}(x)+\sum\limits _{m\in J}H^n_{k-m}P_{m,0}(x)
-\sum\limits _{m=-M}^{M}G^{(k,n)}_mP_{m,S(m)+1}(x) ,
\end{equation}
where 
$$\tilde{P}_{k,n}(x)=H^n_{k-x}\chi _{(-\infty,-M]\cap [M,+\infty )}(x)+
T_2+H^n_{k-x}\chi _{J_0}(x)$$
(the third term is due to the second sum in (\ref{c01})).
Using (\ref{T2}), we can combine first two terms to obtain

\vskip 1mm
{\bf for $n>0$:}
\begin{equation}   \label{Papprox}
\tilde{P}_{k,n}(x)=D^{(k,n)}_{x,S(x)+1}
+D^{(k,n)}_{x-1,S(x)}\delta _{S(x),S(x-1)+1}+D^{(k,n)}_{x+1,S(x)}
\delta _{S(x),S(x+1)+1}+H^n_{k-x}\chi _{J_0}(x).
\end{equation}

\vskip 1mm
  
In order to obtain analogous results for the case $n<0$, we remind that 
expression (\ref{cn3}) contains two terms: $\gamma ^{(n)}_l$
and $\gamma ^{(-n)}_{-l}$. In the previous treatment we have used
only the first term. It means that analogous results for $n<0$ can be
easily obtained by ``reflection'' of all ``ordinates'' with respect to 
horizontal axis, 

\vskip 2mm
{\bf for $n<0$:}
\begin{equation}   \label{prob2*}
P_{k,n}(x)=\tilde{P}^*_{k,n}(x)+\sum\limits _{m\in J}H^{-n}_{k-m}P_{m,0}(x)
-\sum\limits _{m=-M}^{M}G^{(k,n)*}_mP_{m,S(m)+1}(x) ,
\end{equation}
with
\begin{equation}     \label{Papprox2*}
\tilde{P}^*_{k,n}(x)=\chi _{(-M,M)}(x)\biggl(D^{(k,-n)}_{x,-S(x)-1}
+D^{(k,-n)}_{x-1,-S(x)}\delta _{S(x),S(x-1)-1}+
\end{equation}
$$D^{(k,-n)}_{x+1,-S(x)}\delta _{S(x),S(x+1)-1}+H^{-n}_{k-x}\chi _{J_0}(x) \biggr)  ,$$
$$G^{(k,n)*}_m=D^{(k,-n)}_{m,-S(m)}+\delta _{S(m),S(m+1)+1}D^{(k,-n)}_{m+1,-S(m+1)}
+\delta _{S(m),S(m-1)+1}D^{(k,-n)}_{m-1,-S(m-1)} .$$
In (\ref{Papprox2*}) there appears function $\chi _{(-M,M)}(x)$,
because for $n<0$ there is no contribution $P_{k,n}^{(0)}(x)$
due to plane ``tails'' (see remark at the end of Section 3.4).

Note that we cannot represent analogous expressions (\ref{Papprox}) and
(\ref{Papprox2*}) uniquely by writing $|n|$ and $|S(x)|$. 
It is due to the fact that functions $P_{m,n}$ in the upper half plane 
($n>0$) have no influence on functions $P_{m,n}$ in the lower half plane ($n<0$) 
(except through the ground functions), and vice versa. For example, in the sum of 
near-boundary functions (the last term in (\ref{prob2}) and (\ref{prob2*})) 
coefficients $G^{(k,n)}_m$ should be equal to $0$ if $n>0$ and $S(x)\leq 0$ 
or if $n<0$ and $S(x)\geq 0$.

\vskip 1mm

$\tilde{P}_{k,n}(x)$ can be considered as first approximation to 
$P_{k,n}(x)$. Note that a priori there is no reason to neglect the 
second and the third terms in (\ref{prob2}). Normally, we should take
these terms into account, thus the relation (\ref{prob2}) is considered
as a system of linear equations on the near-boundary and 
ground functions. In order to find these functions, we should
close the system of linear equation. Note that 
equations (\ref{prob2}) must be satisfied for any $k$, $n$ 
and $x$, and we can choose appropriate values of $k$ and $n$.
To close the system for near-boundary functions, we take
$\{ (k,n)\: : \: k\in [-M,M],\: n=S(k)+1\}$, i.e. for any
$k\in [-M,M]$
\begin{equation}    \label{eqn}
P_{k,S(k)+1}(x)=\tilde{P}_{k,S(k)+1}(x)+
\sum\limits _{m\in J}H^{(S(k)+1)}_{k-m}P_{m,0}(x)
-\sum\limits _{m=-M}^{M}G^{(k,S(k)+1)}_mP_{m,S(m)+1}(x).
\end{equation}
To close the system for ground functions, we can choose
different conditions. For example, if we consider surfaces with $S(x)>0$
on $x\in (-M,M)$, there are no ground functions, thus there is
no additional condition other than (\ref{eqn}). For the general
case (where $J\ne \emptyset $), we propose the following condition
\begin{equation}   \label{ground}
P_{k,0}(x)=\frac14\biggl(P_{k+1,0}(x)+P_{k-1,0}(x)+P_{k,1}(x)+P_{k,-1}(x)\biggr),
\hskip 5mm k\in J.
\end{equation}
We started from this relation for all the external points (see (\ref{Pidentity})).
Here we just demand that this relation remains valid if we substitute
our approximations for $P_{k,1}(x)$ and $P_{k,-1}(x)$.

\vskip 10mm
\section{ Some numerical verifications }

In this section we briefly present some numerical results 
to check the validity of the approximation.

First of all, in Fig.1 we depict the function $\varphi (\theta )$ 
which plays a central role in this work.

\begin{center}
\includegraphics[width=8cm]{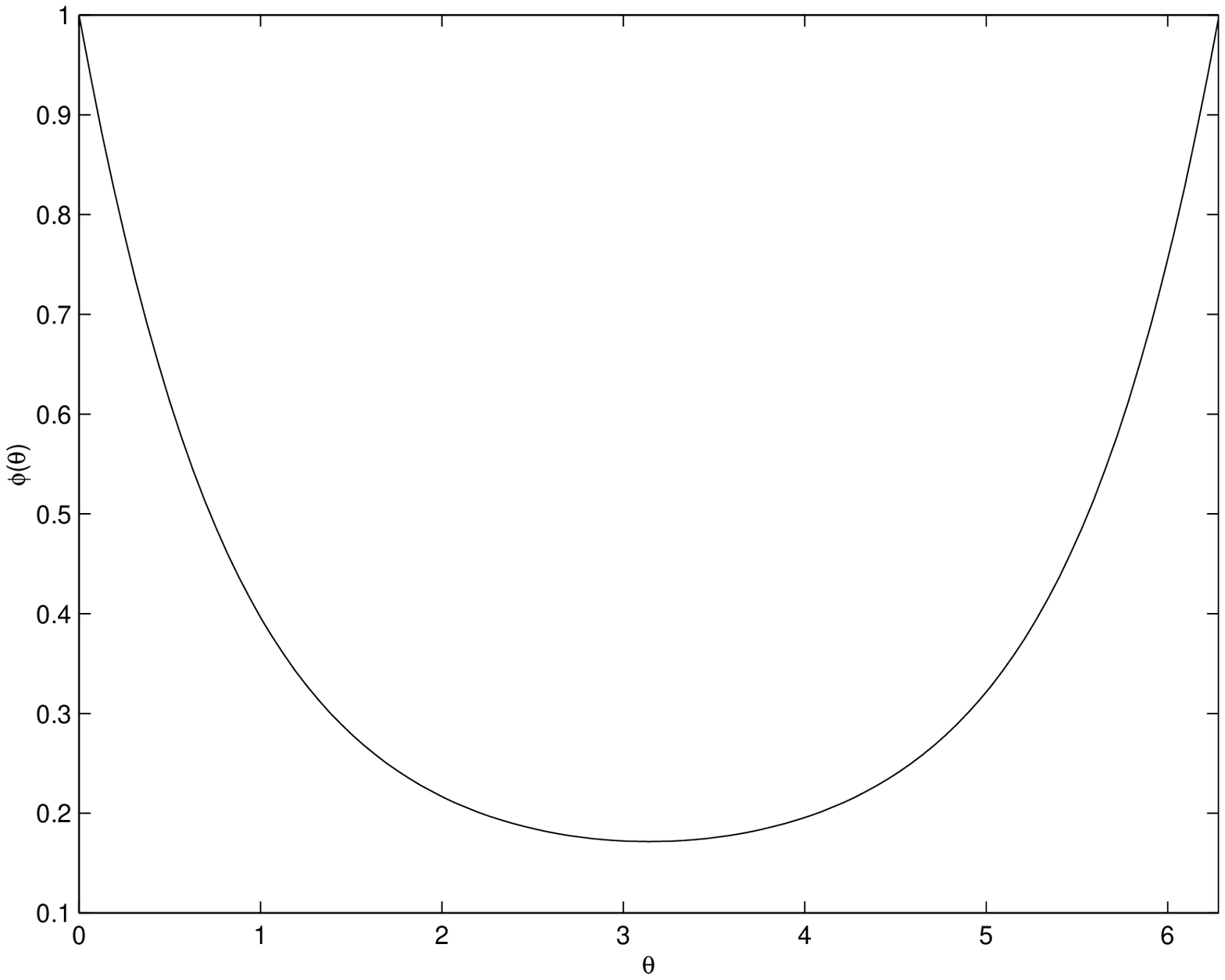}
\end{center}
{\em Figure 1. Function $\varphi (\theta )$. It changes in range 
from $\varphi (0)=1$ to $\varphi (\pi )=3-\sqrt{8}$.}
\vskip 2mm

We can see that $\varphi (\theta )\ll 1$ if $\theta $ is not
in vicinity of $2\pi m$ ($m\in Z$). This property
was used in the limits (\ref{approx}) and (\ref{approx*}).
%
%

Let us consider again the planar surface, $S(x)=0$. Without
simulations, we easily obtain that
$$P_{k,n}(x)=\tilde{P}_{k,n}(x)=H^n_{k-x} .$$
So, in this trivial case our approximation gives the exact result
(cf. (\ref{H})).

\vskip 1mm 

For a particular non-trivial surface the accuracy of the formula 
(\ref{prob2}) can be obtained by comparing its values with numerical 
simulations of random walks. We have taken a simple surface 
represented in Fig.2.

\begin{center}
\includegraphics[width=7cm,height=4cm]{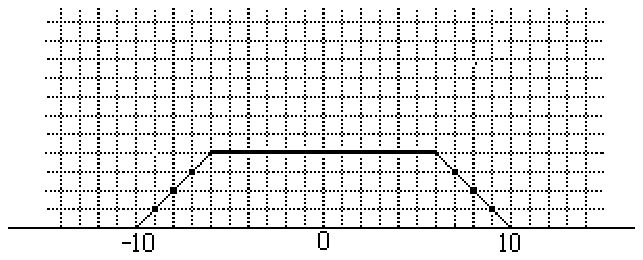}
\end{center}
{\em Figure 2. A simple surface with : $N=4$, $\N =0$, $M=10$.}
\vskip 2mm

In this case there are no ground functions. On the contrary,
there are $21$ near-boundary functions which can be calculated
with the help of (\ref{eqn}). We present two distributions of hitting 
probabilities obtained numerically and through formula (\ref{prob2})
(see Fig.3).

\vskip 1mm

We can conclude that our approximation is quite good.

\begin{center}
\includegraphics[width=12cm,height=10cm]{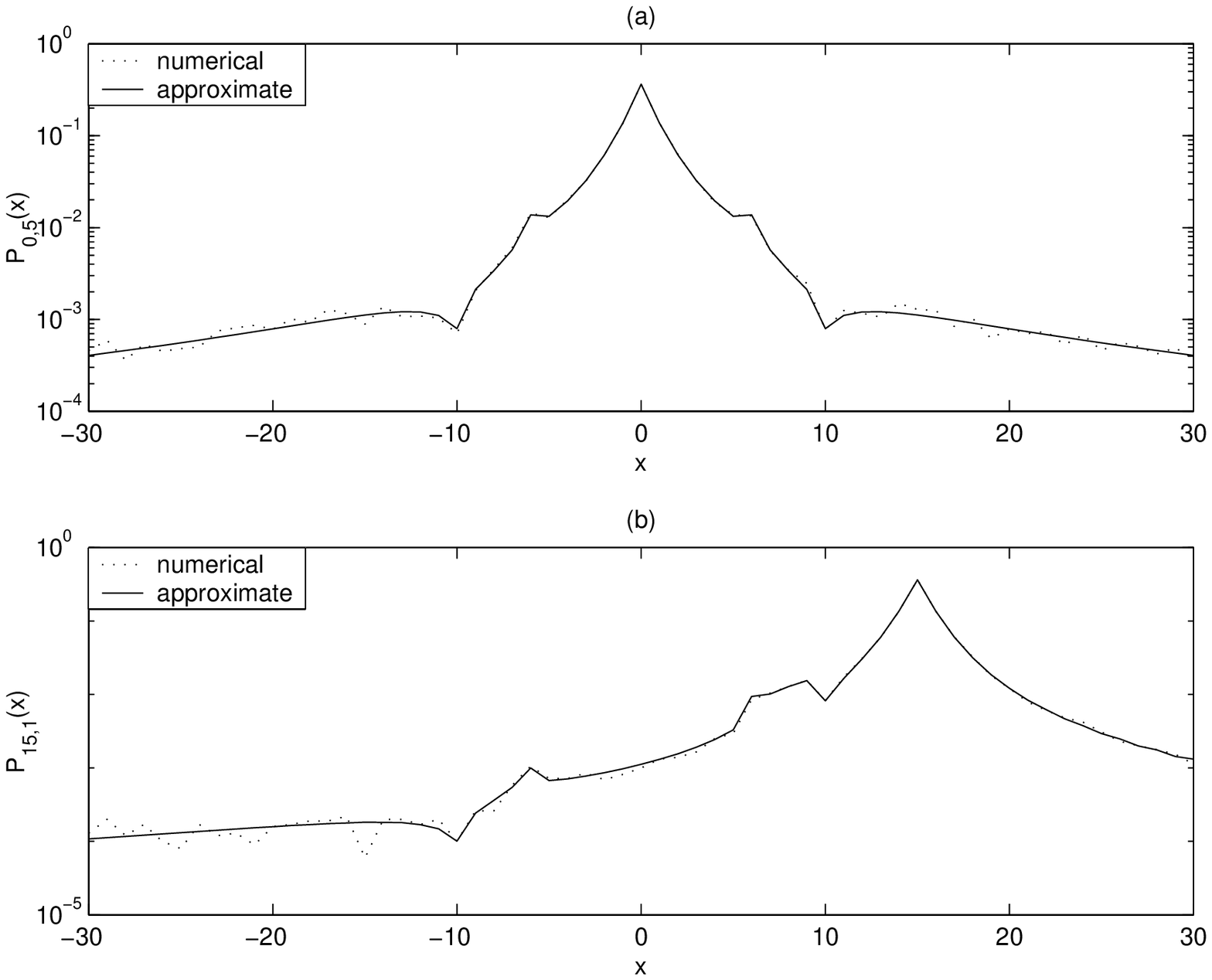}
\end{center}
{\em Figure 3. Distributions of hitting probabilities (in log-scale): 
the probability of the first contact with point $(x,S(x))$
of the surface if started from the point $(0,5)$ (Fig. 3a) or $(15,1)$
(Fig. 3b). }

\vskip 10mm
\section{ Conclusions and possible generalizations }

Let us sum up what has been done. Using the same technique as for a
planar surface, we obtain the reccurence relations for coefficients
$c_n$. For the lower half plane we impose $c_{-N_l}=0$
for a sufficiently large $N_l$. This condition tells that random
walks cannot penetrate through the surface. For the upper half
plane there is no such condition. However, for a sufficiently
large $N_u$ we can use an {\it approximate} condition 
$c_{N_u+1}=\varphi c_{N_u}$ supposing that from a remoted point
the regular surface with compact support looks like a translationally
invariant object. Then we find the explicit solution for 
reccurence relations under these conditions. The influence of
our approximation becomes more and more negligible with
increasing of $N_u$. Taking the limits $N_u\to \infty $ and
$N_l\to \infty $, we express $c_n$ in terms of explicit functions 
$\varphi $ and $\gamma ^{(n)}_l$ and coefficients $\{ \Delta c_l\} $.
Note that two limits $f^{(\infty )}_n$ and $\tilde{f}^{(\infty )}_n$
for the upper and lower half planes are identical. It means that
in the upper half plane we could use the condition $c_{N_u+1}=0$ 
instead of $c_{N_u+1}=\varphi c_{N_u}$. In other words, one could
imagine an absorbing line $y=N_u+1$, and then send it to infinity
(limit $N_u\to \infty $). These two possible approximations
give the same final result (\ref{approx}).
  
Having made these calculations, we express $\Delta c_n$ in terms of
near-boundary functions and combinations of exponential functions.
Finally, we obtain a system of linear equations (\ref{eqn}) 
for near-boundary functions $P_{m,S(m)+1}(x)$ and ground functions $P_{m,0}(x)$.
It can be solved, and after that one can use the approximation
(\ref{prob2}) for any point $(k,n)$.

Numerical analysis shows that this approximation works quite good.

The main conclusion is that we have found an approximate distribution
of hitting probabilities for a rather general surface, pending certain
conditions. In particular, one can make use of these results for a 
further study of the Laplacian transfer problems.

Now one needs to study the role of conditions
which were imposed in the first section. As we said above, 
the compactness condition is the most important. It tells that

-- the surface has a compact support, i.e. there is only {\it finite}
``perturbation'' of the planar surface ;

-- the plane ``tails'' have the same height (zero of the vectical axis).

If we want to consider a surface with {\it infinite} support, we
can obtain the same results but with an infinity 
of near-boundary functions. Thus the system (\ref{eqn}) has 
infinitely many equations, and we cannot proceed any further.
The same difficulty appears if the plane ``tails'' have 
different height: while we step down from the $N_u$-th level to
the level of a lower ``tail'', we must pass through the level of
a higher ``tail''. It means that there appears again an infinity of 
near-boundary functions.
Only if we step down from the $N_u$-th level to the zeroth level 
(and from the $(-N_l)$-th level to the zeroth level), we can avoid 
the appearance of an infinity of near-boundary functions.

The regularity condition is used to simplify certain expressions. 
Nevertheless, it does not seem to be restrictive.
Normally, to apply the technique of characteristic functions, we enumerate 
all sites (points) of the surface. To simplify 
the problem, one can make one of the following assumptions:

-- either suppose that the surface obeys the regularity condition;

-- or be interested in the total hitting probability $P^{(total)}_{k,n}(x)$ 
of points $(x,S(x))$, $(x,S(x)-1)$, ..., $(x,S(x-1)+1)$ (if we authorize
changes of $S(x)$ by more than one unit). In other words, we could
identify the points of the surface which have the same $x$-coordinate.
Any of these assumptions allows to enumerate the points of the surface
with their $x$ coordinate using function $S(x)$.
In the first assumption we consider the surfaces having only 
one point for each $x$; in the second assumption the surfaces 
can have some points with the same $x$, but we are interested in the 
total probability for each $x$.

In order to generalize the method, one can introduce another
{\it parametrization} of the surface. One possible generalization
will be presented in our forthcoming paper.

\section{  Appendices }
\subsection{ Coefficients $H^n_k$ } 

As we have seen, coefficients $H^n_k$ play a central role in
all calculations of hitting probabilities in 2D. Here we briefly 
present some useful properties of $H^n_k$. 
In real form (\ref{H}) becomes,
\begin{equation}   \label{H1}
H^n_k=\int\limits _0^{\pi }\frac{d\theta }{\pi }\cos (k\theta )\varphi ^n(\theta ) .
\end{equation}
We can write two inequalities for $\theta \in (0,\pi )$ :
$$e^{-\theta }\leq \varphi (\theta )\leq e^{-\theta }\frac{1}{\cos (\theta /2)} ,$$
which can be useful for estimations.

Now we are going to calculate the asymptotics of $H^n_k$ for large $k$.
Integrating the expression (\ref{H1}) by parts four times and using 
the values of the derivatives $\varphi ^{k}(s)$ at the points $0$ and $\pi $
(see Table 1), we obtain the asymptotic behaviour 
%
\begin{equation}    \label{asymp}
H^n_k=\frac{n}{\pi k^2}-\frac{n(n^2-0.5)}{\pi k^4}+O(k^{-6}), \hskip 5mm k\gg 1 .
\end{equation}
\begin{center}
\begin{tabular}{| c | c | c | c | c |}  \hline
 $\theta $  & $\varphi $ & $\varphi '$ & $\varphi ''$ & $\varphi '''$  \\ \hline
 $0$        &   $1$     &   $-1$     &   $1$       &      $-1/2$     \\
$\pi $      & $3-\sqrt{8}$&  $0$     &  $3\sqrt{2}/4-1$ &  $0$       \\ \hline
\end{tabular}
\end{center}
{\em Table 1. The values of the derivatives $\varphi ^{(k)}(s)$
at the points $0$ and $\pi $.}

The formula (\ref{asymp}) works rather well for $k\geq 10$. 
\begin{center}
\begin{tabular}{| c | c | c | c | c | c | c |}  \hline
 $n\backslash k$   &  $0$  &  $1$  &  $2$  &  $3$  &  $4$  &  $5$ \\ \hline
 $1$   & $0.3633$ & $0.1366$ & $0.0609$ & $0.0319$ & $0.0189$ & $0.0124$ \\
 $2$   & $0.1803$ & $0.1221$ & $0.0756$ & $0.0477$ & $0.0315$ & $0.0219$ \\
 $3$   & $0.1136$ & $0.0958$ & $0.0715$ & $0.0517$ & $0.0376$ & $0.0278$ \\
 $4$   & $0.0826$ & $0.0759$ & $0.0631$ & $0.0501$ & $0.0392$ & $0.0307$ \\
 $5$   & $0.0651$ & $0.0620$ & $0.0548$ & $0.0464$ & $0.0384$ & $0.0315$ \\ \hline
\end{tabular}
\end{center}
{\em Table 2. The values of coefficients $H^n_k$ for small $k$
($n$ in range from $1$ to $5$).}
\vskip 2mm

Table 2 shows values of $H^n_k$ for small $k$.
The asymptotics of $H^n_k$ for large $n$ is
$$H^n_k=\frac{n}{\pi (n^2+k^2)}+O(n^{-3}) ,$$
i.e. we obtained the same behaviour as for the brownian motion. 
It is quite a reasonable result : if we look on the surface from 
a remoted point, there is no difference between continuous and 
discrete cases.

\vskip 5mm
\subsection{ Manipulation with coefficients $\alpha _n$ and $\beta _n$ }

Here we present some properties of coefficients $\alpha _n$
and $\beta _n$. Also we prove the formula (\ref{cn}).
Using only the definition (\ref{alpha}), we find
\begin{equation}    \label{property}
\alpha _n=\alpha _l\alpha _{n-l+1}-\alpha _{l-1}\alpha _{n-l}
\end{equation}
for any $l\leq n$. Also using (\ref{beta}), we have
$$\beta _n=\alpha _{n+1}-\varphi \alpha _n .$$

\vskip 1mm

Let us prove (\ref{cn}). According to (\ref{solu}), we have
$$c_n=\frac{\beta _{N_u-n}}{\beta _{N_u}}\left(c_0+\sum\limits _{l=1}^n\alpha _l\Delta c_l\right)
+\frac{1}{\beta _{N_u}}\sum\limits _{l=1}^{N_u-n}\Delta c_{n+l}(\beta _{N_u-n}\alpha _{n+l}
-\beta _{N_u}\alpha _l) .$$
Now we should simplify the difference in brackets in the last sum.
\begin{equation}   \label{difference}
\beta _{N_u-n}\alpha _{n+l}-\beta _{N_u}\alpha _l=(\alpha _{N_u-n+1}-
\varphi \alpha _{N_u-n})\alpha _{n+l}-(\alpha _{N_u+1}-\varphi \alpha _{N_u})\alpha _l .
\end{equation}
Consider the difference $\Delta =\alpha _{N_u-n+1}\alpha _{n+l}-\alpha _{N_u+1}\alpha _l$.
Using the property (\ref{property}), we can reduce the index $(n+l)$ in the
first term and $(N_u+1)$ -- in the second term,
$$\Delta =\alpha _{N_u-n+1}(\alpha _l\alpha _{n+1}-\alpha _{l-1}\alpha _n)-
(\alpha _{n+1}\alpha _{N_u-n+1}-\alpha _n\alpha _{N_u-n})\alpha _l=$$
$$=\alpha _n(\alpha _{N_u-n}\alpha _l-\alpha _{N_u-n+1}\alpha _{l-1})=
\alpha _n\alpha _{N_u-n-l+1} $$
(we used the property (\ref{property}) in the last equality).
Thus, we can represent (\ref{difference}) as
$$\beta _{N_u-n}\alpha _{n+l}-\beta _{N_u}\alpha _l=\alpha _n\alpha _{N_u-n-l+1}-
\varphi \alpha _n\alpha _{N_u-n-l}=\alpha _n\beta _{N_u-n-l},$$
hence we find the formula (\ref{cn}),
$$c_n=f^{(N_u)}_n\left(c_0+\sum\limits _{l=1}^n\alpha _l\Delta c_l\right)+\alpha _n
\sum\limits _{l=n+1}^{N_u}f^{(N_u)}_l\Delta c_l .$$


\end{document}